\documentclass[12pt,reqno]{amsart}
\usepackage{enumerate, latexsym, amsmath, amsfonts, amssymb, amsthm, color}
\def\pmod #1{\ ({\rm{mod}}\ #1)}
\def\Z{\Bbb Z}

\def\l{\left}
\def\r{\right}
\def\bg{\bigg}
\def\({\bg(}
\def\){\bg)}
\def\t{\text}
\def\f{\frac}
\def\mo{{\rm{mod}\ }}

\def\ls{\leqslant}

\def\bi{\binom}
\def\al{\alpha}

\def\ve{\varepsilon}

\def\eq{\equiv}

\def\FF#1#2#3{{}_2F_1\bigg(\bmatrix{#1}\\{#2}\end{bmatrix}\bigg|#3\bigg)}
\def\Proof{\noindent{\it Proof}}

\theoremstyle{plain}
\newtheorem{theorem}{Theorem}

\newtheorem{corollary}{Corollary}
\newtheorem{conjecture}{Conjecture}
\theoremstyle{definition}

\theoremstyle{remark}
\newtheorem{remark}{Remark}

 \vspace{4mm}

\begin{document}

\hbox{Nanjing Univ. J. Math. Biquarterly 32(2015), no.\,2, 189--218.}
\medskip

\title
[{New series for some special values of $L$-functions}]
{New series for some special values of $L$-functions}

\author
[Zhi-Wei Sun] {Zhi-Wei Sun}

\address {Department of Mathematics, Nanjing
University, Nanjing 210093, People's Republic of China}
\email{zwsun@nju.edu.cn}

\thanks{
\newline \indent Supported by the National Natural Science Foundation (grant 11171140)
 of China.}

 \begin{abstract} Dirichlet's $L$-functions are natural extensions of the Riemann zeta function. In this paper we first give a brief survey of
 Ap\'ery-like series for some special values of the zeta function and certain $L$-functions. Then, we establish two theorems
 on transformations of certain kinds of congruences. Motivated by the results and based on our computation, we pose 48 new conjectural series
 (most of which involve harmonic numbers) for such special values and related constants. For example, we conjecture that
 \begin{align*}\sum_{k=1}^\infty\f1{k^4\bi{2k}k}\(\f1k+\sum_{j=k}^{2k}\f1j\)=&\f{11}9\zeta(5),
 \\\sum_{k=1}^\infty\f{(-1)^{k-1}}{k^3\bi{2k}k}\(\f1{5k^3}+\sum_{j=1}^{k}\f1{j^3}\)=&\f{2}5\zeta(3)^2,
 \end{align*}
 and
 $$\sum_{k=1}^\infty\f{48^k}{k(2k-1)\bi{4k}{2k}\bi{2k}k}=\f{15}2\sum_{k=1}^\infty\f{(\f k3)}{k^2},$$
 where $(\f k3)$ denotes the Legendre symbol.
\end{abstract}

\maketitle

\section{Introduction}
\setcounter{lemma}{0}
\setcounter{theorem}{0}
\setcounter{corollary}{0}
\setcounter{remark}{0}
\setcounter{equation}{0}
\setcounter{conjecture}{0}

The Riemann zeta function given by
$$\zeta(s)=\sum_{n=1}^\infty\f1{n^s}\ \quad\t{for}\ \Re(s)>1$$
plays an important role in analytic number theory. As Euler proved,
$$2\zeta(2m)=(-1)^{m-1}\f{(2\pi)^{2m}}{(2m)!}B_{2m}\quad\t{for}\ m=1,2,3,\ldots$$
(cf. \cite[pp.\,111-112]{EL}),
where $B_0,B_1,B_2,\ldots$ are Bernoulli numbers defined by
$$ B_0=1,\ \t{and}\ \sum_{k=0}^n\bi{n+1}kB_k=0\ \t{for}\ n=1,2,3,\ldots.$$
In particular,
$$\zeta(2)=\f{\pi^2}6,\ \zeta(4)=\f{\pi^4}{90} \ \t{and}\ \zeta(6)=\f{\pi^6}{945}.$$
It is well known that $\pi$ is transcendental (cf. \cite[pp.\,173-176]{EL}) and thus all those $\zeta(2m)\ (m=1,2,3,\ldots)$
are irrational. In 1978 R. Ap\'ery (cf. \cite{A} and \cite{vP}) successfully established the irrationality of $\zeta(3)$
by using the series
\begin{equation}\label{1.1}\sum_{k=1}^\infty\f{(-1)^{k-1}}{k^3\bi{2k}k}=\f25\zeta(3)
\end{equation}
which converges at a geometric rate with ratio $1/4$ since
$$\bi{2k}k\sim\f{4^k}{\sqrt{k\pi}}\quad\t{as}\ k\to+\infty.$$
In fact, the last identity was first deduced by A. A. Markov \cite{M} in 1890.
There are also some fast converging series for $\zeta(2)$ similar to (\ref{1.1}), e.g.,
\begin{equation}\label{1.2}\sum_{k=1}^\infty\f1{k^2\bi{2k}k}=\f{\pi^2}{18},\ \sum_{k=1}^\infty\f{2^k}{k^2\bi{2k}k}=\f{\pi^2}8,
\ \sum_{k=1}^\infty\f{3^k}{k^2\bi{2k}k}=\f{2}9\pi^2
\end{equation}
(see R. Matsumoto \cite{Ma}). More generally,
\begin{equation}\label{1.3}\arcsin^2 \l(\f x2\r)=\f12\sum_{k=1}^\infty\f{x^{2k}}{k^2\bi{2k}k}\quad\ \t{for}\ |x|\ls2
\end{equation}
(see, e.g., \cite{BC}). It is also well known that
\begin{equation}\label{1.4}\sum_{k=1}^\infty\f1{k^4\bi{2k}k}=\f{17}{36}\zeta(4)=\f{17}{3240}\pi^4
\end{equation}
(cf. \cite[p.\,89]{C}).

Recall that the {\it harmonic numbers} are given by
$$H_n:=\sum_{0<k\ls n}\f1k\ \ (n=0,1,2,\ldots).$$
For $m=2,3,4,\ldots$, we call those numbers
$$H_n^{(m)}:=\sum_{0<k\ls n}\f1{k^m}\ \ (n=0,1,2,\ldots)$$
{\it harmonic numbers of order $m$}. In contrast with (\ref{1.3}), it is also known that for every $n=2,3,\ldots$ we have
$$\f{\arcsin^{2n}(x/2)}{(2n)!}=\f1{4^n}\sum_{k=1}^\infty\f{x^{2k}}{k^2\bi{2k}k}\sum_{0<k_1<k_2<\cdots<k_{n-1}<k}\f{1}{(k_1k_2\cdots k_{n-1})^2}$$
for $|x|\ls 2$, and in particular
\begin{equation}\label{1.5}\arcsin^4\f x2=\f32\sum_{k=1}^\infty\f{H_{k-1}^{(2)}}{k^2\bi{2k}k}x^{2k}\quad\t{for}\ |x|\ls2.
\end{equation}
(See \cite{BC}.) Putting $x=1$ in (\ref{1.5}) we get
\begin{equation*}\sum_{k=1}^\infty\f{H_{k-1}^{(2)}}{k^2\bi{2k}k}=\f{\pi^4}{1944}.
\end{equation*}
Taking derivatives of both sides of the identity in (\ref{1.5}), we obtain
\begin{equation}\label{1.6}\f{4}{\sqrt{4-x^2}}\arcsin^3\f x2=3\sum_{k=1}^\infty\f{H^{(2)}_{k-1}}{k\bi{2k}k}x^{2k-1}\quad\t{for}\ |x|<2.
\end{equation}
It is well known that
$$\arcsin\f x2=\f12\sum_{k=0}^\infty\f{\bi{2k}k}{(2k+1)16^k}x^{2k+1}\quad\t{for}\ |x|\ls2;$$
in particular,
\begin{equation}\label{1.7}\sum_{k=0}^\infty\f{\bi{2k}k}{(2k+1)16^k}=\f{\pi}3\ \ \t{and}\ \ \sum_{k=0}^\infty\f{\bi{2k}k}{(2k+1)8^k}=\f{\sqrt2}4\pi.
\end{equation}
By \cite{BC}, for each $n=1,2,\ldots$ we also have
$$\f{\arcsin^{2n+1}(x/2)}{(2n+1)!}=\f12\sum_{k=0}^\infty\f{\bi{2k}kx^{2k+1}}{(2k+1)16^k}\sum_{0\ls k_1<k_2<\cdots<k_n<k}\f{1}{\prod_{i=1}^n(2k_i+1)^2}$$
provided $|x|\ls 2$; in particular,
$$\arcsin^3\f x2=3\sum_{k=0}^\infty\f{\bi{2k}kx^{2k+1}}{(2k+1)16^k}\sum_{0\ls j<k}\f1{(2j+1)^2}\quad\t{for}\ |x|\ls2.$$
Note that
$$\sum_{0\ls j<k}\f1{(2j+1)^2}=H_{2k}^{(2)}-\f{H_k^{(2)}}4.$$
It is also known that
\begin{equation}\label{1.8}\sum_{k=0}^\infty\f{\bi{2k}k}{(2k+1)^2(-16)^k}=\f{\pi^2}{10}
\end{equation}
(cf. \cite{Ma}) and
\begin{equation}\label{1.9}21\sum_{k=0}^\infty\f{\bi{2k}k}{(2k+1)16^k}\sum_{0\ls j<k}\f1{(2j+1)^2}=\sum_{k=0}^\infty\f{\bi{2k}k}{(2k+1)^316^k}=\f{7\pi^3}{216}
\end{equation}
(cf. I. J. Zucker \cite[(2.23)]{Z} and \cite{HP10}).

M. Koecher \cite{K} and D. Leshchiner \cite[4a]{L} independently proved that for each $n=2,3,\ldots$ we have
\begin{align*}\zeta(2n+1)=&(-1)^{n-1}\sum_{k=1}^\infty\f{(-1)^{k-1}}{k^3\bi{2k}k}
\times\bigg(\f 52\sum_{0<k_1<\cdots<k_{n-1}<k}\f1{\prod_{0<j<n}k_j^2}
\\&+2\sum_{0<m<n}\f{(-1)^{m}}{k^{2m}}\sum_{0<k_1<\cdots<k_{n-1-m}<k}\f1{\prod_{0<j<n-m}k_j^2}\bigg);
\end{align*}
in particular,
\begin{equation}\label{1.10}\zeta(5)=\sum_{k=1}^\infty\f{(-1)^{k-1}}{k^3\bi{2k}k}\l(\f 2{k^2}-\f 52H_{k-1}^{(2)}\r).
\end{equation}
Leshchiner \cite{L} deduced this via a sophisticated analytic method as well as an elegant combinatorial approach.
Inspired by this result, J. Borwein and D. Bradley \cite{BB} used the PSLQ algorithm for finding integer relations to discover that
\begin{equation}\label{1.11}\zeta(7)=\f52\sum_{k=1}^\infty\f{(-1)^{k-1}}{k^3\bi{2k}k}\l(5H^{(4)}_{k-1}+\f1{k^4}\r).
\end{equation}
An extension of this was proved by G. Almkvist and A. Granville \cite{AG} in 1999.
Leshchiner \cite[4b]{L} also proved that for any $n=1,2,3,\ldots$ we have
\begin{align*}&(-1)^{n-1}(1-4^{-n})\zeta(2n)
\\=&\sum_{k=0}^\infty\f{\bi{2k}k}{(2k+1)^2(-16)^k}\bigg(\f54\sum_{1\ls k_1<\cdots<k_{n-1}\ls k}\f1{\prod_{0<j<n}(2k_j-1)^2}
\\&+\sum_{0<m<n}\f{(-1)^m}{(2k+1)^{2m}}\sum_{1\ls k_1<\cdots<k_{n-1-m}\ls k}\f1{\prod_{0<j<n-m}(2k_j-1)^2}\bigg);
\end{align*}
where the last sum is regarded as $1$ when $m=n-1$; in particular,
$$\sum_{k=0}^\infty\f{\bi{2k}k}{(2k+1)^2(-16)^k}\(5\sum_{0\ls j<k}\f1{(2j+1)^2}-\f 4{(2k+1)^2}\)=-\f{\pi^4}{24}.$$

Let $\chi$ be a Dirichlet character modulo a positive integer $m$. The Dirichlet $L$-function associated with the character $\chi$ is given by
$$L(s,\chi):=\sum_{k=1}^\infty\f{\chi(k)}{k^s}\quad\t{for}\ \Re(s)>1.$$
The Dirichlet beta function is defined by
$$\beta(s)=L\l(s,\l(\f{-4}{\cdot}\r)\r)=\sum_{k=0}^\infty\f{(-1)^k}{(2k+1)^s}\quad\t{for}\ \Re(s)>0,$$
where $(-)$ denotes the Kronecker symbol. As Euler observed,
$$\beta(2n+1)=\f{(-1)^nE_{2n}}{4^{n+1}(2n)!}\pi^{2n+1}\quad\t{for all}\ n=0,1,2,\ldots$$
(cf. (3.63) of \cite[p.\,112]{EL}), where $E_0,E_1,E_2,\ldots$ are Euler numbers defined by
$$E_0=1,\ \t{and}\ \sum^n_{k=0\atop 2\mid k}\bi nkE_{n-k}=0\ \ \t{for}\ n=1,2,3,\ldots.$$
In particular,
$$\sum_{k=0}^\infty\f{(-1)^k}{2k+1}=\f{\pi}4,\ \sum_{k=0}^\infty\f{(-1)^k}{(2k+1)^3}=\f{\pi^3}{32},\ \sum_{k=0}^\infty\f{(-1)^k}{(2k+1)^5}=\f{5\pi^5}{1536}.$$
S. Ramanujan found that for $0<|x|\ls \pi/2$ we have
$$\sum_{k=0}^\infty\f{\bi{2k}k\sin^{2k+1}x}{(2k+1)^24^k}=x\log|2\sin x|+\sum_{k=1}^\infty\f{\sin(2kx)}{2k^2}$$
(cf. B. C. Berndt and P. T. Joshi \cite[p.\, 89]{BJ}), consequently
$$\sum_{k=0}^\infty\f{\bi{2k}k}{(2k+1)^28^k}=\f{\pi}{4\sqrt2}\log2+\f{G}{\sqrt2}$$
and
$$\sum_{k=0}^\infty\f{\bi{2k}k}{(2k+1)^216^k}=\f{3\sqrt3}{4}K$$
(cf. \cite[pp.\,39-40]{BJ}), where
$$G:=\beta(2)=\sum_{k=0}^\infty\f{(-1)^k}{(2k+1)^2}$$
is the Catalan constant, and
$$K:=L\l(2,\l(\f{-3}{\cdot}\r)\r)=\sum_{k=1}^\infty\f{(\f k3)}{k^2}=0.781302412896\ldots$$
with $(\f k3)$ the Legendre symbol.
In contrast with (\ref{1.1}) and (\ref{1.4}), in 1985 I. J. Zucker \cite{Z} proved the following remarkable identities:
\begin{align*}\sum_{k=1}^\infty\f1{k^3\bi{2k}k}=&\f{\sqrt3}2\pi K-\f 43\zeta(3),
\\\sum_{k=1}^\infty\f{2^k}{k^3\bi{2k}k}=&\f{\pi^2}8\log 2+\pi G-\f{35}{16}\zeta(3),
\\\sum_{k=1}^\infty\f{3^k}{k^3\bi{2k}k}=&\f 29\pi^2\log 3+\f{2\sqrt3}3\pi K-\f{26}9\zeta(3),
\\\sum_{k=0}^\infty\f{\bi{2k}k}{(2k+1)^416^k}=&\f{\pi}{12}\zeta(3)+\f{27\sqrt3}{64}L,
\\\sum_{k=1}^\infty\f1{k^5\bi{2k}k}=&\f{9\sqrt3}8\pi L+\f{\pi^2}9\zeta(3)-\f{19}3\zeta(5),
\end{align*}
where
$$L:=L\l(4,\l(\f{-3}{\cdot}\r)\r)=\sum_{k=1}^\infty\f{(\f k3)}{k^4}.$$
Z.-W. Sun \cite{S11} conjectured that
$$\sum_{k=1}^\infty\f{(15k-4)(-27)^{k-1}}{k^3\bi{2k}k\bi{3k}k}=K
\ \t{and}\ \sum_{k=1}^\infty\f{(5k-1)(-144)^k}{k^3\bi{2k}k^2\bi{4k}{2k}}=-\f{45}2K;$$
the first formula here was confirmed in 2012 by Kh. Hessami Pilehrood and T. Hessami Pilehrood \cite{HP12b}, and the second one
was later proved by J. Guillera and M. Rogers \cite{GR}.

L. van Hamme \cite{vH} investigated corresponding $p$-adic congruences for certain hypergeometric series
involving the Gamma function or $\pi=\Gamma(1/2)^2$. This stimulated later studies of $p$-adic congruences for some well-known series.
For example, R. Tauraso \cite{T} proved that (\ref{1.1}) and the first identity in (\ref{1.2}) have the following $p$-adic analogs:
\begin{equation}\label{1.12}\sum_{k=1}^{p-1}\f1{k^2\bi{2k}k}\eq\f{H_{p-1}}{3p}\pmod{p^3}\ \t{and}\ \sum_{k=1}^{p-1}\f{(-1)^{k-1}}{k^3\bi{2k}k}\eq \f25\cdot\f{H_{p-1}}{p^2}\pmod{p^3},
\end{equation}
where $p>5$ is a prime. (It is known that $H_{p-1}\eq-p^2B_{p-3}/3\pmod{p^3}$ for any prime $p>3$.)
Z.-W. Sun \cite{S11,S14} showed that
$$\sum_{k=1}^{(p-1)/2}\f1{k^2\bi{2k}k}\eq\l(\f{-1}p\r)\f43E_{p-3}
\ \t{and}\ \sum_{k=1}^{(p-1)/2}\f{(-1)^{k-1}}{k^3\bi{2k}k}\eq2B_{p-3}\pmod p$$
for any prime $p>5$.
By \cite[Lemma 2.1]{S11}, if $p$ is an odd prime then
\begin{equation}\label{1.13}k\bi{2k}k\bi{2(p-k)}{p-k}\eq(-1)^{\lfloor 2k/p\rfloor-1}2p\ \pmod{p^2}
\end{equation}
for all $k=1,\ldots,p-1$. (The paper \cite{T} contains a similar technique.)  Motivated by this and (\ref{1.4}), Z.-W. Sun \cite{S11} conjectured that
\begin{align}\label{1.14}p\sum_{k=1}^{p-1}\f1{k^4\bi{2k}k}\eq&\f{H_{p-1}}{p^2}-\f7{45}p^2B_{p-5}\pmod{p^3},
\\\label{1.15}\sum_{k=1}^{p-1}\f{\bi{2k}k}{k^3}\eq&-2\f{H_{p-1}}{p^2}-\f{13}{27}H_{p-1}^{(3)}\pmod{p^4},
\end{align}
where $p>7$ is a prime.
This remains open but their mod $p$ versions have been confirmed
by Kh. Hessami Pilehrood and T. Hessami Pilehrood \cite{HP12a}.
Similarly, motivated by (\ref{1.10}), (\ref{1.11}) and (\ref{1.13}), we conjecture that
\begin{equation}\label{1.16}
\sum_{k=1}^{(p-1)/2}\f{(-1)^k}{k^3\bi{2k}k}\l(H_{k-1}^{(2)}-\f4{5k^2}\r)\eq10B_{p-5}\pmod{p}
\end{equation}
for every prime $p>3$, and that
\begin{equation}\label{1.17}\sum_{k=1}^{(p-1)/2}\f{(-1)^k}{k^3\bi{2k}k}\l(5H_{k-1}^{(4)}+\f1{k^4}\r)\eq58B_{p-7}\pmod p
\end{equation}
and
\begin{equation}\label{1.18}\sum_{k=1}^{p-1}\f{(-1)^k}{k^3\bi{2k}k}\l(5H_{k-1}^{(4)}+\f1{k^4}\r)\eq\f2{35}B_{p-7}\pmod p
\end{equation}
for any prime $p>7$.
Note that if $p$ is an odd prime and $m$ is a positive integer with $(p-1)\nmid m$ then for every $k=1,\ldots,p-1$ we have
$$H_{p-k}^{(m)}=H_{p-1}^{(m)}-\sum_{0<j<k}\f1{(p-j)^m}\eq (-1)^{m-1}H_{k-1}^{(m)}\ \pmod p.$$

Let $p>3$ be a prime. Motivated by the first identity in (\ref{1.7}), Z.-W. Sun \cite{S11j} showed that
\begin{equation}\label{1.19}\sum_{k=0}^{(p-3)/2}\f{\bi{2k}k}{(2k+1)16^k}\eq0\pmod{p^2},\ \sum_{k=(p+1)/2}^{p-1}\f{\bi{2k}k}{(2k+1)16^k}\eq\f p3E_{p-3}\pmod{p^2}.
\end{equation}
In view of (\ref{1.8}) and the last equality in (\ref{1.9}), Z.-W. Sun \cite{S11j} conjectured that if $p>5$ then
\begin{align}\label{1.20}\sum_{k=0}^{(p-3)/2}\f{\bi{2k}k}{(2k+1)^2(-16)^k}\eq&\f{H_{p-1}}{5p}\pmod{p^3},
\\\label{1.21}\sum_{k=(p+1)/2}^{p-1}\f{\bi{2k}k}{(2k+1)^2(-16)^k}\eq&-\f p4B_{p-3}\pmod{p^2},
\\\label{1.22}\sum_{k=0}^{(p-3)/2}\f{\bi{2k}k}{(2k+1)^316^k}\eq&\l(\f{-1}p\r)\l(\f{H_{p-1}}{4p^2}+\f{p^2}{36}B_{p-5}\r)\pmod{p^3}.
\end{align}
The congruences (\ref{1.20}) and (\ref{1.21}) were confirmed in \cite{HPT} and \cite{S14} respectively.

Here we pose a curious conjectural series for $K$ as well as its related $p$-adic congruences.

\begin{conjecture}\label{Conj1.1} {\rm (i)} We have
\begin{equation}\label{1.23}\sum_{k=1}^\infty\f{48^k}{k(2k-1)\bi{4k}{2k}\bi{2k}k}=\f{15}2K.
\end{equation}

{\rm (ii)} For any prime $p>3$, we have
\begin{equation}\label{1.24}\sum_{k=1}^{p-1}\f{\bi{4k}{2k+1}\bi{2k}k}{48^k}\eq\f 5{12}p^2B_{p-2}\l(\f13\r)\pmod{p^3}\end{equation}
and
\begin{equation}\label{1.25}p^2\sum_{k=1}^{p-1}\f{48^k}{k(2k-1)\bi{4k}{2k}\bi{2k}k}\eq 4\l(\f p3\r)+4p\pmod{p^2},
\end{equation}
where $B_{n}(x)$ denotes the Bernoulli polynomial of degree $n$.
\end{conjecture}

\begin{remark}\label{Rem1.1} In \cite[Conj. 5.14]{S11} the author made the following conjecture for any prime $p>3$:
If $p\eq1\ (\mo\ 3)$ and $p=x^2+3y^2$ with $x,y\in\Z$ and $x\eq1\ (\mo\ 3)$, then
 \begin{gather*}\sum_{k=0}^{p-1}\f{\bi{2k}k\bi{4k}{2k}}{48^k}\eq 2x-\f p{2x}\ (\mo\ p^2),
 \\\sum_{k=0}^{p-1}\f{k+1}{48^k}\bi{2k}k\bi{4k}{2k}\eq x\ (\mo\ p^2).
 \end{gather*}
If $p\eq2\ (\mo\ 3)$, then
$$\sum_{k=0}^{p-1}\f{\bi{2k}k\bi{4k}{2k}}{48^k}\eq\f{3p}{2\bi{(p+1)/2}{(p+1)/6}}\pmod{p^2}.$$
\end{remark}

In view of (\ref{1.24}), it is interesting to investigate what primes $p>3$ satisfy the congruence $B_{p-2}(1/3)\eq0\pmod p$.
The first such a prime is $205129$ and the second one must be greater than $2\times10^7$.
By the standard heuristic arguments (cf. \cite[pp.\,28-29]{CP}), it seems that there should be infinitely many such primes.

In the next section we establish two theorems on transformations of certain kinds of congruences.
Motivated by the results and based on our computation via the PSLQ algorithm, we pose in Sections 3-5 many conjectural series involving harmonic numbers or higher-order harmonic numbers
for some special values of $L$-functions,
most of which are mainly motivated by their $p$-adic analogs we found first. All the conjectural series in this paper converge
at geometric rates and so they can be easily checked via some mathematical softwares like {\tt Mathematica}.

\section{Transformations of certain kinds of congruences}
\setcounter{lemma}{0}
\setcounter{theorem}{0}
\setcounter{corollary}{0}
\setcounter{remark}{0}
\setcounter{equation}{0}
\setcounter{conjecture}{0}

\begin{theorem}\label{Th2.1} Let $p$ be an odd prime and let $a$ be an integer with $p\nmid a$. Let $b,c\in\Z$, $m,r\in\Z^+=\{1,2,3,\ldots\}$ and $\ve\in\{\pm1\}$. Then
\begin{equation}\label{2.1}\begin{aligned}&(-1)^r\sum_{k=1}^{p-1}\f{\bi{2k}k}{a^kk^r}\(b\sum_{j=1}^k\f{\ve^j}{j^m}+c\f{\ve^k}{k^m}\)
\\\eq&2bp\sum_{j=1}^{p-1}\f{\ve^j}{j^m}\sum_{k=1}^{p-1}\f{a^{k-1}}{k^{r+1}\bi{2k}k}
\\&-(-1)^{m}2\ve p\sum_{k=1}^{p-1}\f{a^{k-1}}{k^{r+1}\bi{2k}k}\(b\sum_{j=1}^k\f{\ve^j}{j^m}-(b+c)\f{\ve^k}{k^m}\)\pmod p,
\end{aligned}
\end{equation}
and also
\begin{equation}\label{2.2}\begin{aligned}&\f1{(-2)^r}\l(\f{-a}p\r)\sum_{k=1}^{p-1}\f{\bi{2k}k}{a^kk^r}\(b\sum_{j=1}^k\f{\ve^j}{j^m}+c\f{\ve^k}{k^m}\)
\\\eq&b\sum_{j=1}^{(p-1)/2}\f{\ve^j}{j^m}\sum_{k=0}^{(p-3)/2}\f{\bi{2k}k}{(2k+1)^r(16/a)^k}-(-2)^{m}\l(\f{\ve}p\r)
\\&\times\sum_{k=0}^{(p-3)/2}\f{\bi{2k}k}{(2k+1)^r(16/a)^k}\(b\sum_{j=0}^k\f{\ve^j}{(2j+1)^m}-(b+c)\f{\ve^k}{(2k+1)^m}\)\pmod p.
\end{aligned}
\end{equation}
\end{theorem}
\Proof. Let $n=(p-1)/2$. For $k\in\{1,\ldots,p-1\}$, clearly $p\mid\bi{2k}k$ if and only if $k>n$.
By (\ref{1.13}),
$$k\bi{2k}k\eq\f{-2p}{\bi{2(p-k)}{p-k}}\pmod{p}\quad\t{for all}\ k=1,\ldots,n.$$
Thus
\begin{align*}&\sum_{k=1}^{p-1}\f{\bi{2k}k}{a^kk^r}\(b\sum_{j=1}^k\f{\ve^j}{j^m}+c\f{\ve^k}{k^m}\)
\\\eq&\sum_{k=1}^{n}\f{k\bi{2k}k}{a^kk^{r+1}}\(b\sum_{j=1}^k\f{\ve^j}{j^m}+c\f{\ve^k}{k^m}\)
\\\eq&-2p\sum_{k=1}^n\f{1}{a^kk^{r+1}\bi{2(p-k)}{p-k}}\(b\sum_{j=1}^k\f{\ve^j}{j^m}+c\f{\ve^k}{k^m}\)
\\\eq&-2p\sum_{k=1}^{p-1}\f{1}{a^kk^{r+1}\bi{2(p-k)}{p-k}}\(b\sum_{j=1}^k\f{\ve^j}{j^m}+c\f{\ve^k}{k^m}\)
\\\eq&-2p\sum_{k=1}^{p-1}\f{1}{a^{p-k}(p-k)^{r+1}\bi{2k}k}\(b\sum_{j=1}^{p-k}\f{\ve^j}{j^m}+c\f{\ve^{p-k}}{(p-k)^m}\)
\\\eq&-2p\sum_{k=1}^{p-1}\f{a^{k-1}}{(-k)^{r+1}\bi{2k}k}\(b\sum_{j=1}^{p-1}\f{\ve^j}{j^m}-b\sum_{0<s<k}\f{\ve^{p-s}}{(p-s)^m}+c\f{\ve^{k-1}}{(-k)^m}\)
\\\eq&2p\sum_{k=1}^{p-1}\f{(-1)^ra^{k-1}}{k^{r+1}\bi{2k}k}\(b\sum_{j=1}^{p-1}\f{\ve^j}{j^m}-(-1)^m\ve\(b\sum_{s=1}^k\f{\ve^{s}}{s^m}-(b+c)\f{\ve^{k}}{k^m}\)\)
\pmod p.
\end{align*}
This proves (\ref{2.1}).

For each $k=0,\ldots,n$, we have
\begin{equation}\label{2.3}\f{\bi{2k}k}{16^k}\eq\l(\f{-1}p\r)\bi{2(n-k)}{n-k}\pmod p
\end{equation}
since
\begin{align*}\bi{2k}k=\bi{-1/2}k(-4)^k\eq&\bi nk(-4)^k=\bi n{n-k}(-4)^k
\\\eq&\bi{-1/2}{n-k}(-4)^k=\f{\bi{2(n-k)}{n-k}}{(-4)^{n-k}}(-4)^k
\\\eq&(-1)^n\bi{2(n-k)}{n-k}16^k\pmod p.
\end{align*}
Therefore
\begin{align*}&\sum_{k=1}^{p-1}\f{\bi{2k}k}{a^kk^r}\(b\sum_{j=1}^k\f{\ve^j}{j^m}+c\f{\ve^k}{k^m}\)
\\\eq&\sum_{k=1}^{n}\f{\bi{2k}k}{a^kk^r}\(b\sum_{j=1}^k\f{\ve^j}{j^m}+c\f{\ve^k}{k^m}\)
\\=&\sum_{k=0}^{n-1}\f{\bi{2(n-k)}{n-k}}{a^{n-k}(n-k)^{r}}\(b\sum_{j=1}^{n-k}\f{\ve^j}{j^m}+c\f{\ve^{n-k}}{(n-k)^m}\)
\\\eq&\sum_{k=0}^{n-1}\l(\f {-a}p\r)\f{a^k\bi{2k}k}{16^k(-1/2-k)^{r}}\(b\sum_{j=1}^n\f{\ve^j}{j^m}-b\sum_{0\ls s<k}\f{\ve^{n-s}}{(n-s)^m}+c\f{\ve^{n-k}}{(-1/2-k)^m}\)
\\\eq&\l(\f{-a}p\r)(-2)^r\sum_{k=0}^{n-1}\f{\bi{2k}k}{(2k+1)^r(16/a)^k}
\\&\times\(b\sum_{j=1}^n\f{\ve^j}{j^m}-\ve^n(-2)^m\(b\sum_{0\ls s<k}\f{\ve^s}{(2s+1)^m}-c\f{\ve^k}{(2k+1)^m}\)\)
\\\eq&\l(\f{-a}p\r)(-2)^r\sum_{k=0}^{n-1}\f{\bi{2k}k}{(2k+1)^r(16/a)^k}
\\&\times\(b\sum_{j=1}^n\f{\ve^j}{j^m}-(-2)^m\l(\f{\ve}p\r)\(b\sum_{s=0}^k\f{\ve^s}{(2s+1)^m}-(b+c)\f{\ve^k}{(2k+1)^m}\)\)
\pmod p.
\end{align*}
This proves (\ref{2.2}).

In view of the above, we have completed the proof of Theorem \ref{Th2.1}. \qed

\begin{remark}\label{Rem2.1} Let $p$ be an odd prime, and let $m\in\Z^+$ and $\ve\in\{\pm1\}$. It is well known that $H_{p-1}^{(m)}\eq0\pmod p$ when $(p-1)\nmid m$.
Note also that
$$\sum_{k=1}^{p-1}\f{\ve^k}{k^m}=\sum_{k=1}^{(p-1)/2}\(\f{\ve^k}{k^m}+\f{\ve^{p-k}}{(p-k)^m}\)\eq(1+(-1)^m\ve)\sum_{k=1}^{(p-1)/2}\f{\ve^k}{k^m}\pmod p$$
and hence $\sum_{k=1}^{p-1}(-1)^k/k^m\eq0\pmod p$ if $2\mid m$.
\end{remark}

\begin{corollary}\label{Cor2.1} Let $p>5$ be a prime. For any $b,c\in\Z$, $m\in\Z^+$ and $\ve\in\{\pm1\}$, we have
\begin{equation}\label{2.4}\begin{aligned}&\sum_{k=1}^{p-1}\f{\bi{2k}k}{k}\(b\sum_{j=1}^k\f{\ve^j}{j^m}+c\f{\ve^k}{k^m}\)
\\\eq&(-1)^{m}2\ve p\sum_{k=1}^{p-1}\f{1}{k^{2}\bi{2k}k}\(b\sum_{j=1}^k\f{\ve^j}{j^m}-(b+c)\f{\ve^k}{k^m}\)\pmod p,
\end{aligned}
\end{equation}
and also
\begin{equation}\label{2.5}\begin{aligned}&\f{(-1)^m}{2^{m+1}}\l(\f{-\ve}p\r)\sum_{k=1}^{p-1}\f{\bi{2k}k}{k}\(b\sum_{j=1}^k\f{\ve^j}{j^m}+c\f{\ve^k}{k^m}\)
\\\eq&\sum_{k=0}^{(p-3)/2}\f{\bi{2k}k}{(2k+1)16^k}\(b\sum_{j=0}^k\f{\ve^j}{(2j+1)^m}-(b+c)\f{\ve^k}{(2k+1)^m}\)\pmod p.
\end{aligned}
\end{equation}
\end{corollary}
\Proof. In view of (\ref{1.12}) and (\ref{1.19}),
$$p\sum_{k=1}^{p-1}\f1{k^2\bi{2k}k}\eq\sum_{k=0}^{(p-3)/2}\f{\bi{2k}k}{(2k+1)16^k}\eq0\pmod p.$$
Thus, by applying Theorem \ref{Th2.1} with $a=r=1$ we immediately obtain (\ref{2.4}) and (\ref{2.5}). \qed

\begin{corollary}\label{Cor2.2} Let $p>5$ be a prime. For any $b,c\in\Z$, $m\in\Z^+$ and $\ve\in\{\pm1\}$, we have
\begin{equation}\label{2.6}\begin{aligned}&\sum_{k=1}^{p-1}\f{(-1)^{k-1}\bi{2k}k}{k^2}\(b\sum_{j=1}^k\f{\ve^j}{j^m}+c\f{\ve^k}{k^m}\)
\\\eq&(-1)^{m}2\ve p\sum_{k=1}^{p-1}\f{(-1)^{k-1}}{k^3\bi{2k}k}\(b\sum_{j=1}^k\f{\ve^j}{j^m}-(b+c)\f{\ve^k}{k^m}\)\pmod p,
\end{aligned}
\end{equation}
and also
\begin{equation}\label{2.7}\begin{aligned}&\f{(-1)^m}{2^{m+2}}\l(\f{\ve}p\r)\sum_{k=1}^{p-1}\f{(-1)^{k-1}\bi{2k}k}{k^2}\(b\sum_{j=1}^k\f{\ve^j}{j^m}+c\f{\ve^k}{k^m}\)
\\\eq&\sum_{k=0}^{(p-3)/2}\f{\bi{2k}k}{(2k+1)^2(-16)^k}\(b\sum_{j=0}^k\f{\ve^j}{(2j+1)^m}-(b+c)\f{\ve^k}{(2k+1)^m}\)\pmod p.
\end{aligned}
\end{equation}
\end{corollary}
\Proof. In view of (\ref{1.12}) and (\ref{1.20}),
$$p\sum_{k=1}^{p-1}\f{(-1)^{k-1}}{k^3\bi{2k}k}\eq\sum_{k=0}^{(p-3)/2}\f{\bi{2k}k}{(2k+1)^2(-16)^k}\eq0\pmod p.$$
Applying Theorem \ref{Th2.1} with $a=-1$ and $r=2$, we immediately obtain (\ref{2.6}) and (\ref{2.7}). \qed

\begin{corollary}\label{Cor2.3} Let $p>3$ be a prime. For any $b,c\in\Z$, $m\in\Z^+$ and $\ve\in\{\pm1\}$, we have
\begin{equation}\label{2.8}\begin{aligned}&\sum_{k=1}^{p-1}\f{\bi{2k}k}{k^3}\(b\sum_{j=1}^k\f{\ve^j}{j^m}+c\f{\ve^k}{k^m}\)-\f 23bB_{p-3}\sum_{k=1}^{p-1}\f{\ve^k}{k^m}
\\\eq&(-1)^{m}2\ve p\sum_{k=1}^{p-1}\f{1}{k^4\bi{2k}k}\(b\sum_{j=1}^k\f{\ve^j}{j^m}-(b+c)\f{\ve^k}{k^m}\)\pmod p,
\end{aligned}
\end{equation}
and also
\begin{equation}\label{2.9}\begin{aligned}&\f{(-1)^m}{2^{m+3}}\l(\f{-\ve}p\r)\(\sum_{k=1}^{p-1}\f{\bi{2k}k}{k^3}\(b\sum_{j=1}^k\f{\ve^j}{j^m}+c\f{\ve^k}{k^m}\)-\f23bB_{p-3}\sum_{k=1}^{(p-1)/2}\f{\ve^k}{k^m}\)
\\\eq&\sum_{k=0}^{(p-3)/2}\f{\bi{2k}k}{(2k+1)^316^k}\(b\sum_{j=0}^k\f{\ve^j}{(2j+1)^m}-(b+c)\f{\ve^k}{(2k+1)^m}\)\pmod p.
\end{aligned}
\end{equation}
\end{corollary}
\Proof. Let $n=(p-1)/2$. By (\ref{2.3}),
\begin{align*}\sum_{k=0}^{n-1}\f{\bi{2k}k}{(2k+1)^316^k}\eq&\l(\f{-1}p\r)\sum_{k=0}^{n-1}\f{\bi{2(n-k)}{n-k}}{(2k+1)^3}=\l(\f{-1}p\r)\sum_{k=1}^n\f{\bi{2k}k}{(2(n-k)+1)^3}
\\\eq&\l(\f{-1}p\r)\sum_{k=1}^n\f{\bi{2k}k}{(-2k)^3}\eq-\f18\l(\f{-1}p\r)\sum_{k=1}^{p-1}\f{\bi{2k}k}{k^3}\pmod p.
\end{align*}
As conjectured by the author and proved in \cite{HP12a},
$$p\sum_{k=1}^{p-1}\f1{k^4\bi{2k}k}\eq\f{H_{p-1}}{p^2}\eq-\f{B_{p-3}}3\pmod p$$
and $$\sum_{k=1}^{p-1}\f{\bi{2k}k}{k^3}\eq-2\f{H_{p-1}}{p^2}\eq\f23B_{p-3}\pmod p.$$
So
\begin{equation}\label{2.10}\sum_{k=0}^{(p-3)/2}\f{\bi{2k}k}{(2k+1)^316^k}\eq-\l(\f{-1}p\r)\f{B_{p-3}}{12}\pmod p,
\end{equation}
which confirms the mod $p$ version of (\ref{1.22}).
Applying Theorem \ref{Th2.1} with $a=1$ and $r=3$, we immediately obtain (\ref{2.8}) and (\ref{2.9}). \qed

\begin{theorem}\label{Th2.2} Let $p$ be an odd prime and let $m$ be an integer with $p\nmid m(m-4)$.
Let $\al$ be a positive integer, and let $a_0,a_1,\ldots,a_{p^{\al}-1}$ be $p$-adic integers. Define $a_k^*=\sum_{j=0}^k\bi kj(-1)^ja_j$ for $k=0,1,\ldots,p^{\al}-1$.
Then we have the congruences
\begin{equation}\label{2.11}\sum_{k=0}^{p^{\al}-1}\f{\bi{2k}k}{(4-m)^k}a_k^*\eq\l(\f{m(m-4)}{p^{\al}}\r)\sum_{k=0}^{p^{\al}-1}\f{\bi{2k}k}{m^k}a_k\pmod{p}
\end{equation}
and
\begin{equation}\label{2.12}m\sum_{k=0}^{p^{\al}-1}\f{k\bi{2k}k}{(4-m)^k}a_k^*\eq\l(\f{m(m-4)}{p^{\al}}\r)\sum_{k=0}^{p^{\al}-1}((m-4)k-2)\f{\bi{2k}k}{m^k}a_k\pmod{p},
\end{equation}
where $(\f{\cdot}{p^{\al}})$ denotes the Jacobi symbol.
\end{theorem}
\Proof. Let $n=(p^{\al}-1)/2$. By Lucas' theorem (see, e.g., \cite{HS}), for each $k=n+1,\ldots,p^a-1$ we have
$$\bi{2k}k=\bi{p^a+(2k-p^a)}{0\cdot p^a+k}\eq\bi{1}0\bi{2k-p^a}k=0\pmod p.$$
Also, for any $j,k\in\{0,\ldots,n\}$ with $j+k\ls n$, we have
\begin{align*}\f{\bi{n-j}k}{\bi{-1/2-j}k}=&\prod_{0\ls i<k}\f{(p^\al-1)/2-j-i}{-1/2-j-i}
\\=&\prod_{0\ls i<k}\l(1-\f{p^{\al}}{2(i+j)+1}\r)\eq1\pmod{p}.
\end{align*}

As $\bi{-1/2}k=\bi{2k}k/(-4)^k$ for $k=0,1,2,\ldots$, we have
\begin{align*}&\sum_{k=0}^n\f{\bi{2k}k}{(4-m)^k}\sum_{j=0}^k\bi kj(-1)^ja_j
\\=&\sum_{j=0}^n(-1)^ja_j\sum_{k=j}^n\bi{-1/2}k\l(\f{-4}{4-m}\r)^k\bi kj
\\=&\sum_{j=0}^n(-1)^ja_j\bi{-1/2}j\sum_{k=j}^n\bi{-1/2-j}{k-j}\l(\f 4{m-4}\r)^k
\\=&\sum_{j=0}^n\f{\bi{2j}j}{4^j}a_j\l(\f4{m-4}\r)^j\sum_{i=0}^{n-j}\bi{-1/2-j}i\l(\f 4{m-4}\r)^i
\end{align*}
and hence
\begin{align*}\sum_{k=0}^{p^{\al}-1}\f{\bi{2k}k}{(4-m)^k}a_k^*
\eq&\sum_{j=0}^n\f{\bi{2j}j}{(m-4)^j}a_j\sum_{i=0}^{n-j}\bi{n-j}i\l(\f 4{m-4}\r)^i
\\=&\sum_{j=0}^n\f{\bi{2j}j}{(m-4)^j}a_j\l(1+\f4{m-4}\r)^{n-j}
\\=&\sum_{j=0}^n\f{\bi{2j}j}{(m-4)^j}a_j\l(\f m{m-4}\r)^n\l(\f{m-4}m\r)^j
\\=&\sum_{j=0}^n\f{\bi{2j}j}{m^j}a_j\l(\f m{m-4}\r)^{\f{p-1}2\sum_{r=0}^{\al-1}p^r}
\\\eq&\l(\f{m(m-4)}{p^{\al}}\r)\sum_{j=0}^{p^{\al}-1}\f{\bi{2j}j}{m^j}a_j\pmod p.
\end{align*}
This proves (\ref{2.11}). Similarly,
\begin{align*}&\sum_{k=0}^n\f{k\bi{2k}k}{(4-m)^k}\sum_{j=0}^k\bi kj(-1)^ja_j
\\=&\sum_{j=0}^n(-1)^ja_j\bi{-1/2}j\sum_{k=j}^nk\bi{-1/2-j}{k-j}\l(\f 4{m-4}\r)^k
\\=&\sum_{j=0}^n\f{\bi{2j}j}{4^j}a_j\l(\f4{m-4}\r)^j\sum_{i=0}^{n-j}(j+i)\bi{-1/2-j}i\l(\f 4{m-4}\r)^i
\end{align*}
and hence
\begin{align*}\sum_{k=0}^{p^{\al}-1}\f{k\bi{2k}k}{(4-m)^k}a_k^*
\eq&\sum_{j=0}^n\f{\bi{2j}ja_j}{(m-4)^j}\sum_{i=0}^{n-j}(j+i)\bi{n-j}i\l(\f 4{m-4}\r)^i
\\=&\sum_{j=0}^n\f{j\bi{2j}ja_j}{(m-4)^j}\l(1+\f4{m-4}\r)^{n-j}
\\&+\sum_{j=0}^n\f{(n-j)\bi{2j}j}{(m-4)^j}a_j\sum_{0<i\ls n-j}\bi{n-j-1}{i-1}\l(\f 4{m-4}\r)^i
\\=&\sum_{j=0}^n\f{j\bi{2j}j}{(m-4)^j}a_j\l(\f m{m-4}\r)^{n-j}
\\&+\sum_{j=0}^n\f{(n-j)\bi{2j}j}{(m-4)^j}a_j\l(1+\f4{m-4}\r)^{n-j-1}\f 4{m-4}
\\=&\sum_{j=0}^n\f{j\bi{2j}ja_j}{(m-4)^j}\l(\f{m-4}m\r)^{j-n}
\\&+\f4m\sum_{j=0}^n\f{(n-j)\bi{2j}j}{(m-4)^j}a_j\l(\f{m-4}m\r)^{j-n}
\\\eq&\l(\f{m(m-4)}{p^{\al}}\r)\sum_{j=0}^{p^{\al}-1}\f{\bi{2j}j}{m^j}a_j\l(j+\f 4m\l(-\f12-j\r)\r)\pmod p.
\end{align*}
So (\ref{2.12}) also holds. This concludes the proof. \qed

\begin{remark}\label{Rem2.2} Theorem \ref{Th2.2} was motivated by the identity (1.18) in \cite{S14a}.
\end{remark}

\begin{corollary}\label{Cor2.4} For $k=0,1,2,\ldots$ define
$$f_k(x):=\sum_{j=0}^k\bi kj^2\bi{2j}kx^j\quad\ \t{and}\ \quad g_k(x):=\sum_{j=0}^k\bi kj^2\bi{2j}jx^j.$$
Let $p$ be an odd prime and let $\al$ be any positive integer. Then, for any integer $m\not\eq0,-4\pmod p$, we have
\begin{equation}\label{2.13}\sum_{k=0}^{p^{\al}-1}\f{\bi{2k}kg_k(x)}{(m+4)^k}\eq\l(\f{m(m+4)}{p^{\al}}\r)\sum_{k=0}^{p^{\al}-1}\f{\bi{2k}k}{m^k}f_k(x)\pmod p
\end{equation}
and
\begin{equation}\label{2.14}
m\sum_{k=0}^{p^{\al}-1}\f{k\bi{2k}kg_k(x)}{(m+4)^k}\eq\l(\f{m(m+4)}{p^{\al}}\r)\sum_{k=0}^{p^{\al}-1}((m+4)k+2)\f{\bi{2k}k}{m^k}f_k(x)\pmod p.
\end{equation}
\end{corollary}
\Proof. By \cite[Theorem 2.2]{S15a},
$$g_k(x)=\sum_{j=0}^k\bi kj(-1)^j((-1)^jf_j(x))\quad\t{for all}\ k=0,1,2,\ldots.$$
Thus, by applying Theorem \ref{Th2.2} we obtain that
$$\sum_{k=0}^{p^{\al}-1}\f{\bi{2k}kg_k(x)}{(4-(-m))^k}\eq\l(\f{-m(-m-4)}{p^{\al}}\r)\sum_{k=0}^{p^{\al}-1}\f{\bi{2k}k}{(-m)^k}(-1)^kf_k(x)\pmod p$$
and that
\begin{align*}&-m\sum_{k=0}^{p^{\al}-1}\f{k\bi{2k}k}{(4-(-m))^k}g_k(x)
\\\eq&\l(\f{-m(-m-4)}{p^{\al}}\r)\sum_{k=0}^{p^{\al}-1}((-m-4)k-2)\f{\bi{2k}k}{(-m)^k}(-1)^kf_k(x)\pmod p.
\end{align*}
So, both (\ref{2.13}) and (\ref{2.14}) hold. \qed

\begin{remark}\label{Rem2.3} The author \cite[(2.6)]{S13} proved that
$$\sum_{k=0}^{p-1}\f{\bi{2k}kf_k(x)}{(-4)^k}\eq\sum_{k=0}^{p-1}\f{\bi{2k}k^3}{16^k}x^k\pmod{p^2}\quad\t{for any prime}\ p>3.$$
\end{remark}

\section{Conjectural formulas involving ordinary harmonic numbers}
\setcounter{lemma}{0}
\setcounter{theorem}{0}
\setcounter{corollary}{0}
\setcounter{remark}{0}
\setcounter{equation}{0}
\setcounter{conjecture}{0}

\begin{conjecture}\label{Conj3.1} {\rm (i)} We have
\begin{align}\label{3.1}\sum_{k=1}^\infty\f{H_{2k}+2/(3k)}{k^2\bi{2k}k}=&\zeta(3),
\\\label{3.2}\sum_{k=1}^\infty\f{H_{2k}+2H_k}{k^2\bi{2k}k}=&\f53\zeta(3),
\\\label{3.3}\sum_{k=1}^\infty\f{H_{2k}+17H_k}{k^2\bi{2k}k}=&\f52\sqrt3\ \pi K.
\end{align}

{\rm (ii)} Let $p>3$ be a prime. Then
\begin{align*}\sum_{k=1}^{(p-1)/2}\f{3H_{2k}+2/k}{k^2\bi{2k}k}\eq& B_{p-3}\pmod p,
\\\sum_{k=1}^{p-1}\f{3H_{2k}+2/k}{k^2\bi{2k}k}\eq&-4\f{H_{p-1}}{p^2}-\f9{10}p^2B_{p-5}\pmod{p^3},
\\\sum_{k=1}^{p-1}\f{H_{2k}+2H_k}{k^2\bi{2k}k}\eq&-\f 83\cdot\f{H_{p-1}}{p^2}-\f{17}{30}p^2B_{p-5}\pmod{p^3},
\\p\sum_{k=1}^{p-1}\f{H_{2k}+17H_k}{k^2\bi{2k}k}\eq&\f54\l(\f p3\r)B_{p-2}\l(\f13\r)\pmod p.
\end{align*}
\end{conjecture}
\begin{remark}\label{Rem3.1} A combination of (\ref{3.1}) and (\ref{3.2}) yields
\begin{equation*}\sum_{k=1}^\infty\f{3H_k-1/k}{k^2\bi{2k}k}=\zeta(3)
\end{equation*}
for which {\tt Mathematica 9} could yield a ``{\tt proof}" after running the {\tt FullSimplify} command half an hour.
Combining (\ref{3.1})-(\ref{3.3}) we find exact values of
$$\sum_{k=1}^\infty \f1{k^3\bi{2k}k},\ \sum_{k=1}^\infty\f{H_k}{k^2\bi{2k}k}\ \t{and}\ \sum_{k=1}^\infty\f{H_{2k}}{k^2\bi{2k}k}.$$
Recently G.-S. Mao and Z.-W. Sun \cite{MS} showed that for any prime $p>3$ we have
\begin{align*}\sum_{k=1}^{p-1}\f{\bi{2k}k}kH_k\eq&\f13\l(\f p3\r)B_{p-2}\l(\f 13\r)\pmod p,
\\\sum_{k=1}^{p-1}\f{\bi{2k}k}kH_{2k}\eq&\f7{12}\l(\f p3\r)B_{p-2}\l(\f 13\r)\pmod p.
\end{align*}
\end{remark}

\begin{conjecture}\label{Conj3.2} {\rm (i)} We have
\begin{align}\label{3.4}\sum_{k=1}^\infty\f{2^k}{k^2\bi{2k}k}\l(H_{\lfloor k/2\rfloor}-(-1)^k\f 2k\r)=&\f 74\zeta(3),
\\\label{3.5}\sum_{k=1}^\infty\f{2^k}{k^2\bi{2k}k}\l(2H_{2k}-3H_k+\f2k\r)=&\f 74\zeta(3),
\\\label{3.6}\sum_{k=1}^\infty\f{2^k}{k^2\bi{2k}k}\l(6H_{2k}-11H_k+\f 8k\r)=&2\pi G,
\\\label{3.7}\sum_{k=1}^\infty\f{2^k}{k^2\bi{2k}k}\l(2H_{2k}-7H_k+\f2k\r)=&-\f{\pi^2}2\log2,
\\\label{3.8}\sum_{k=1}^\infty\f{3^k}{k^2\bi{2k}k}\l(6H_{2k}-8H_k+\f5k\r)=&\f{26}3\zeta(3),
\\\label{3.9}\sum_{k=1}^\infty\f{3^k}{k^2\bi{2k}k}\l(6H_{2k}-10H_k+\f7k\r)=&2\sqrt3\,\pi K,
\\\label{3.10}\sum_{k=1}^\infty\f{3^k}{k^2\bi{2k}k}\l(H_k+\f1{2k}\r)=&\f{\pi^2}3\log3.
\end{align}
And
\begin{align}\label{3.11}\sum_{k=1}^\infty\f{L_{2k}}{k^2\bi{2k}k}\l(\f1k+\f1{k+1}+\cdots+\f1{2k}\r)=&\f{41\zeta(3)+4\pi^2\log\phi}{25},
\\\label{3.12}\sum_{k=1}^\infty\f{v_k}{k^2\bi{2k}k}\l(\f1k+\f1{k+1}+\cdots+\f1{2k}\r)=&\f{124\zeta(3)+\pi^2\log\l(5^5\phi^6\r)}{50},
\end{align}
where $\phi$ is the famous golden ratio $(\sqrt5+1)/2$, the Lucas numbers $L_0,L_1,L_2,\ldots$
are given by
$$L_0=2,\ L_1=1,\ \t{and}\ L_{n+1}=L_n+L_{n-1}\ \ \t{for all}\ n=1,2,3,\ldots,$$
and $v_0,v_1,v_2,\ldots,$ are defined by
$$v_0=2,\ v_1=5,\ \t{and}\ v_{n+1}=5(v_n-v_{n-1})\ \ \t{for all}\ n=1,2,3,\ldots.$$

{\rm (ii)} For any prime $p>3$ we have
\begin{align*}p^2\sum_{k=1}^{p-1}\f{2^k}{k^2\bi{2k}k}\l(2H_{2k}-3H_k+\f2k\r)\eq&-2q_p(2)-H_{p-1}\pmod{p^3},
\\p^2\sum_{k=1}^{p-1}\f{2^k}{k^2\bi{2k}k}\l(6H_{2k}-11H_k+\f8k\r)\eq&-6q_p(2)+\l(\f{-1}p\r)2pE_{p-3}\pmod{p^2},
\\p^2\sum_{k=1}^{p-1}\f{2^k}{k^2\bi{2k}k}\l(2H_{2k}-7H_k+\f2k\r)\eq&-2q_p(2)+2p\,q_p(2)^2\pmod{p^2},
\\\sum_{k=1}^{p-1}\f{\bi{2k}k}{k2^k}\l(2H_{2k}-3H_k\r)\eq&\f{49}{24}pB_{p-3}\pmod{p^2},
\\p^2\sum_{k=1}^{p-1}\f{3^k}{k^2\bi{2k}k}\l(6H_{2k}-8H_k+\f5k\r)\eq&-9q_p(3)-6H_{p-1}\pmod{p^3},
\\p^2\sum_{k=1}^{p-1}\f{3^k}{k^2\bi{2k}k}\l(6H_{2k}-10H_k+\f7k\r)\eq&-9q_p(2)+\l(\f{p}3\r)\f p2B_{p-2}\l(\f13\r)\pmod{p^2},
\\\sum_{k=1}^{p-1}\f{\bi{2k}k}{k3^k}\l(3H_{2k}-4H_k\r)\eq&\f{26}{9}pB_{p-3}\pmod{p^2},
\end{align*}
where $q_p(a)$ for $a\not\eq0\pmod p$ denotes the Fermat quotient $(a^{p-1}-1)/p$.
\end{conjecture}
\begin{remark}\label{Rem3.2} (a) Combining (\ref{3.5})-(\ref{3.7}) we find the exact values of
$$\sum_{k=1}^\infty\f{2^k}{k^2\bi{2k}k},\ \ \sum_{k=1}^\infty\f{2^kH_k}{k^2\bi{2k}k}\ \ \t{and}\ \ \sum_{k=1}^\infty\f{2^kH_{2k}}{k^2\bi{2k}k}.$$
Similarly, combining (\ref{3.8})-(\ref{3.10}) we get the exact values of
$$\sum_{k=1}^\infty\f{3^k}{k^2\bi{2k}k},\ \ \sum_{k=1}^\infty\f{3^kH_k}{k^2\bi{2k}k}\ \ \t{and}\ \ \sum_{k=1}^\infty\f{3^kH_{2k}}{k^2\bi{2k}k}.$$
The author \cite[Remark 5.2]{S15} used {\tt Mathematica 7} to find that
$$\sum_{k=1}^\infty\f{4^kH_{k-1}}{k^2\bi{2k}k}=7\zeta(3)\ \ \t{and}\ \ \sum_{k=1}^\infty\f{4^kH_{2k-1}}{k^2\bi{2k}k}=\f{21}2\zeta(3).$$

(b) Applying (\ref{1.3}) with $x=\phi,1-\phi$ we obtain
$$\sum_{k=1}^\infty\f{\phi^{2k}}{k^2\bi{2k}k}=2\arcsin^2\f{\sqrt5+1}4=2\l(\f{3\pi}{10}\r)^2=\f{9\pi^2}{50}$$
and
$$\sum_{k=1}^\infty\f{(1-\phi)^{2k}}{k^2\bi{2k}k}=2\arcsin^2\f{1-\sqrt5}4=2\l(-\f{\pi}{10}\r)^2=\f{\pi^2}{50}.$$
Therefore,
$$\sum_{k=1}^\infty\f{L_{2k}}{k^2\bi{2k}k}=\sum_{k=1}^\infty\f{\phi^{2k}+(1-\phi)^{2k}}{k^2\bi{2k}k}=\f{\pi^2}5.$$
(\ref{1.3}) with $x=\sqrt{(5\pm\sqrt5)/2}$ yields that
$$2\l(\f{2\pi}5\r)^2=\sum_{k=1}^\infty\f{((5+\sqrt5)/2)^k}{k^2\bi{2k}k}\ \t{and}\ \ 2\l(\f{\pi}5\r)^2=\sum_{k=1}^\infty\f{((5-\sqrt5)/2)^k}{k^2\bi{2k}k}.$$
So we have
$$\sum_{k=1}^\infty\f{v_k}{k^2\bi{2k}k}=\sum_{k=1}^\infty\f1{k^2\bi{2k}k}\(\(\f{5+\sqrt5}2\)^k+\(\f{5-\sqrt5}2\)^k\)=\f25\pi^2.$$
\end{remark}

\begin{conjecture}\label{Conj3.3} {\rm (i)} We have
\begin{align}\label{3.13}\sum_{k=1}^\infty(-1)^{k-1}\f{10H_k-3/k}{k^3\bi{2k}k}=&\f{\pi^4}{30},
\\\label{3.14}\sum_{k=1}^\infty(-1)^{k-1}\f{H_{2k}+4H_k}{k^3\bi{2k}k}=&\f2{75}\pi^4,
\\\label{3.15}\sum_{k=1}^\infty\f{H_{2k}-H_k+2/k}{k^4\bi{2k}k}=&\f{11}9\zeta(5),
\\\label{3.16}\sum_{k=1}^\infty\f{3H_{2k}-102H_k+28/k}{k^4\bi{2k}k}=&-\f{55}{18}\pi^2\zeta(3),
\\\label{3.17}\sum_{k=1}^\infty\f{97H_{2k}-163H_k+227/k}{k^4\bi{2k}k}=&\f{165}8\sqrt3\pi L.
\end{align}

{\rm (ii)} Let $p>3$ be a prime. Then
\begin{align*}\sum_{k=1}^{p-1}\f{(-1)^{k-1}}{k^2}\bi{2k}k\l(H_{k-1}+\f3{10k}\r)\eq&\f{16}{25}p^2B_{p-5}\pmod{p^3},
\\\sum_{k=1}^{p-1}\f{(-1)^{k-1}}{k^3\bi{2k}k}\l(10H_k-\f3k\r)\eq&\f25pB_{p-5}\pmod{p^2},
\\\sum_{k=1}^{p-1}\f{(-1)^{k-1}}{k^2}\bi{2k}k\l(H_{2k-1}+4H_{k-1}\r)\eq&\f{184}{25}p^2B_{p-5}\pmod{p^3}\ \t{if}\ p>5,
\\p\sum_{k=1}^{p-1}\f{(-1)^{k-1}(H_{2k}+4H_k)}{k^3\bi{2k}k}\eq&\f25\cdot\f{H_{p-1}}{p^2}+\f 8{25}p^2B_{p-5}\pmod{p^3},
\\\sum_{k=1}^{p-1}\f{\bi{2k}k}{k^3}\l(H_{2k-1}-H_k-\f1k\r)\eq&\f{22}{15}pB_{p-5}\pmod{p^2},
\\p^2\sum_{k=1}^{p-1}\f{H_{2k}-H_k+2/k}{k^4\bi{2k}k}\eq&\f{H_{p-1}}{p^2}+\f49p^2B_{p-5}\pmod{p^3},
\\\sum_{k=1}^{p-1}\f{\bi{2k}k}{k^3}\l(3H_{2k-1}-102H_{k-1}-\f{28}k\r)\eq&\f{44}{15}pB_{p-5}\pmod{p^2},
\\p^2\sum_{k=1}^{p-1}\f{3H_{2k}-102H_k+28/k}{k^4\bi{2k}k}\eq&3\f{H_{p-1}}{p^2}+\f{196}{15}p^2B_{p-5}\pmod{p^3}.
\end{align*}
If $p>5$, then
$$\sum_{k=1}^{p-1}\f{\bi{2k}k}{k^3}\l(97H_{2k-1}-163H_{k-1}-\f{227}k\r)\eq-\f{55}{48}\l(\f p3\r)B_{p-4}\l(\f13\r)\pmod{p}$$
and
$$p^2\sum_{k=1}^{p-1}\f{97H_{2k}-163H_k+227/k}{k^4\bi{2k}k}\eq97\f{H_{p-1}}{p^2}+\f{55}{96}\l(\f p3\r)pB_{p-4}\l(\f13\r)\pmod{p^2}.$$
\end{conjecture}

\begin{conjecture}\label{Conj3.4} {\rm (i)} We have
\begin{align}\label{3.18}\sum_{k=0}^\infty\f{\bi{2k}k}{(2k+1)16^k}\l(3H_{2k+1}+\f4{2k+1}\r)=&8G,
\\\label{3.19}\sum_{k=0}^\infty\f{\bi{2k}k}{(2k+1)^2(-16)^k}\l(5H_{2k+1}+\f{12}{2k+1}\r)=&14\zeta(3),
\\\label{3.20}\sum_{k=0}^\infty\f{\bi{2k}k}{(2k+1)^316^k}\l(9H_{2k+1}+\f{32}{2k+1}\r)=&40\beta(4)+\f5{12}\pi\zeta(3).
\end{align}

{\rm (ii)} Let $p$ be an odd prime. Then
$$\sum_{k=0}^{(p-3)/2}\f{\bi{2k}k}{(2k+1)16^k}\l(3H_{2k+1}+\f4{2k+1}\r)\eq4E_{p-3}\pmod p.$$
If $p>3$, then
\begin{align*}&\sum_{k=0}^{(p-3)/2}\f{\bi{2k}k}{(2k+1)^2(-16)^k}\l(5H_{2k+1}+\f{12}{2k+1}\r)
\\&\qquad\eq-15\f{H_{p-1}}{p^2}-\f{29}{10}p^2B_{p-5}\pmod{p^3}.
\end{align*}
\end{conjecture}

\section{Conjectural formulas involving higher-order harmonic numbers}
\setcounter{lemma}{0}
\setcounter{theorem}{0}
\setcounter{corollary}{0}
\setcounter{remark}{0}
\setcounter{equation}{0}
\setcounter{conjecture}{0}

For any positive integers $k$ and $m$, obviously
$$\sum_{j=1}^k\f{(-1)^j}{j^m}+H_k^{(m)}=\sum_{j=1}^k\f{(-1)^j+1}{j^m}=\sum_{0<j\ls\lfloor k/2\rfloor}\f2{(2j)^m}=2^{1-m}H_{\lfloor k/2\rfloor}^{(m)}.$$

\begin{conjecture}\label{Conj4.1} {\rm (i)} We have
\begin{equation}\label{4.1}\sum_{k=1}^\infty\f{6H_{\lfloor k/2\rfloor}^{(2)}-(-1)^k/k^2}{k^2\bi{2k}k}=\f{13}{1620}\pi^4.
\end{equation}
Also,
\begin{align}\label{4.2}\sum_{k=1}^\infty\f{H_k^{(3)}}{k^2\bi{2k}k}=&\f{\zeta(5)+2\zeta(2)\zeta(3)}9,
\\\label{4.3}\sum_{k=1}^\infty\f{(-1)^k}{k^3\bi{2k}k}\(10\sum_{j=1}^k\f{(-1)^j}{j^2}-\f{(-1)^k}{k^2}\)=&\f{29\zeta(5)-2\zeta(2)\zeta(3)}6,
\\\label{4.4}\sum_{k=1}^\infty\f1{k^2\bi{2k}k}\(24\sum_{j=1}^k\f{(-1)^j}{j^3}-17\f{(-1)^k}{k^3}\)=&7\zeta(5)-6\zeta(2)\zeta(3).
\end{align}

{\rm (ii)} Let $p>3$ be a prime. Then
\begin{gather*}p\sum_{k=1}^{p-1}\f{6H_{\lfloor k/2\rfloor}^{(2)}-(-1)^k/k^2}{k^2\bi{2k}k}\eq22\f{H_{p-1}}{p^2}-\f{97}{45}p^2B_{p-5}\pmod{p^3},
\\\sum_{k=1}^{p-1}\f{H_k^{(3)}}{k^2\bi{2k}k}\eq\f{29}{45}B_{p-5}\pmod p,
\ \sum_{k=1}^{p-1}\f{\bi{2k}k}kH_{k-1}^{(3)}\eq\f2{45}pB_{p-5}\pmod{p^2}.
\end{gather*}
Also,
\begin{align*}\sum_{k=1}^{p-1}\f{(-1)^k}{k^3\bi{2k}k}\(10\sum_{j=1}^k\f{(-1)^j}{j^2}-\f{(-1)^k}{k^2}\)\eq&\f{121}{30}B_{p-5}\pmod p,
\\\sum_{k=1}^{p-1}\f1{k^2\bi{2k}k}\(24\sum_{j=1}^k\f{(-1)^j}{j^3}-17\f{(-1)^k}{k^3}\)\eq&\f{23}5B_{p-5}\pmod p.
\end{align*}
\end{conjecture}

\begin{conjecture}\label{Conj4.2} {\rm (i)} We have
\begin{align}\label{4.5}\sum_{k=1}^\infty(-1)^{k-1}\f{H_k^{(3)}+1/(5k^3)}{k^3\bi{2k}k}=&\f25\zeta(3)^2,
\\\label{4.6}\sum_{k=1}^\infty\f{H_{k-1}^{(2)}-1/k^2}{k^4\bi{2k}k}=&-\f{313\pi^6}{612360},
\\\label{4.7}\sum_{k=1}^\infty\f{3H_k^{(4)}-1/k^4}{k^2\bi{2k}k}=&\f{163\pi^6}{136080}.
\end{align}
Also,
\begin{align}\label{4.8}\sum_{k=1}^\infty\f1{k^2\bi{2k}k}\(8\sum_{j=1}^k\f{(-1)^j}{j^4}+\f{(-1)^k}{k^4}\)=&-\f{97}{34020}\pi^6-\f{22}{15}\zeta(3)^2,
\\\label{4.9}\sum_{k=1}^\infty\f1{k^4\bi{2k}k}\(72\sum_{j=1}^k\f{(-1)^j}{j^2}-\f{(-1)^k}{k^2}\)=&-\f{31}{1134}\pi^6-\f{34}5\zeta(3)^2,
\\\label{4.10}\sum_{k=1}^\infty\f{(-1)^k}{k^3\bi{2k}k}\(40\sum_{0<j<k}\f{(-1)^j}{j^3}-7\f{(-1)^k}{k^3}\)=&-\f{367}{27216}\pi^6+6\zeta(3)^2.
\end{align}

{\rm (ii)} Let $p>3$ be a prime. Then
\begin{align*}p\sum_{k=1}^{p-1}\f{(-1)^{k-1}}{k^3\bi{2k}k}\l(5H_k^{(3)}+\f1{k^3}\r)\eq&2B_{p-5}\pmod p,
\\ p\sum_{k=1}^{p-1}\f{H_{k-1}^{(2)}-1/k^2}{k^4\bi{2k}k}\eq&\f 49B_{p-5}\pmod{p}.
\end{align*}
If $p>7$, then
\begin{align*}
\\\sum_{k=1}^{p-1}\f{\bi{2k}k}{k^5}+\sum_{k=1}^{p-1}\f{\bi{2k}k}{k^3}H_k^{(2)}\eq&-\f{20}{27}\cdot\f{H_{p-1}^{(3)}}{p^2}-\f{197}{756}p^2B_{p-7}\pmod{p^3},
\\\sum_{k=1}^{p-1}\f{\bi{2k}k}k\l(3H_{k-1}^{(4)}+\f1{k^4}\r)\eq&\f56\cdot\f{H_{p-1}^{(3)}}{p^2}+\f{31}{21}p^2B_{p-7}\pmod{p^3},
\\\sum_{k=1}^{p-1}\f{3H_k^{(4)}-1/k^4}{k^2\bi{2k}k}\eq&\f 5{12}\cdot\f{H_{p-1}^{(3)}}{p^2}-\f{499}{336}p^2B_{p-7}\pmod{p^3}.
\end{align*}
When $p>5$, we have
\begin{align*}p\sum_{k=1}^{p-1}\f1{k^2\bi{2k}k}\(\sum_{j=1}^k\f{(-1)^j}{j^4}+\f{(-1)^k}{8k^4}\)\eq&-\f{B_{p-5}}3\pmod p,
\\p\sum_{k=1}^{p-1}\f1{k^4\bi{2k}k}\(72\sum_{j=1}^k\f{(-1)^j}{j^2}-\f{(-1)^k}{k^2}\)\eq&-2B_{p-5}\pmod p,
\\p\sum_{k=1}^{p-1}\f{(-1)^k}{k^3\bi{2k}k}\(40\sum_{0<j<k}\f{(-1)^j}{j^3}-7\f{(-1)^k}{k^3}\)\eq&\f 52B_{p-5}\pmod p.
\end{align*}
\end{conjecture}

\begin{conjecture}\label{Conj4.3}  We have
\begin{align}\label{4.11}\sum_{k=1}^\infty\f{33H_k^{(5)}+4/k^5}{k^2\bi{2k}k}=&-\f{45}8\zeta(7)+\f{13}3\zeta(2)\zeta(5)+\f{85}6\zeta(3)\zeta(4),
\\\label{4.12}\sum_{k=1}^\infty\f{33H_k^{(3)}+8/k^3}{k^4\bi{2k}k}
=&-\f{259}{24}\zeta(7)-\f{98}9\zeta(2)\zeta(5)+\f{697}{18}\zeta(3)\zeta(4),
\end{align}
and
\begin{equation}\label{4.13}\begin{aligned}&\sum_{k=1}^\infty\f{(-1)^k}{k^3\bi{2k}k}\(110\sum_{j=1}^k\f{(-1)^j}{j^4}+29\f{(-1)^k}{k^4}\)
\\=&\f{223}{24}\zeta(7)-\f{301}6\zeta(2)\zeta(5)+\f{221}2\zeta(3)\zeta(4).
\end{aligned}
\end{equation}
\end{conjecture}

\section{Conjectural formulas involving $\sum_{j=0}^k(\pm1)^j/(2j+1)^m$}
\setcounter{lemma}{0}
\setcounter{theorem}{0}
\setcounter{corollary}{0}
\setcounter{remark}{0}
\setcounter{equation}{0}
\setcounter{conjecture}{0}

\begin{conjecture}\label{Conj5.1} {\rm (i)} We have
\begin{align}\label{5.1}\sum_{k=0}^\infty\f{\bi{2k}k}{(2k+1)16^k}\sum_{j=0}^k\f1{(2j+1)^3}=&\f 5{18}\pi\zeta(3),
\\\label{4.2}\sum_{k=0}^\infty\f{\bi{2k}k}{(2k+1)^2(-16)^k}\sum_{j=0}^k\f{(-1)^j}{(2j+1)^2}=&\f{\pi^2G}{10}+\f{\pi\zeta(3)}{240}+\f{27\sqrt3}{640}L.
\end{align}

{\rm (ii)} Let $p>3$ be a prime. Then
\begin{equation*}\sum_{k=0}^{(p-3)/2}\f{\bi{2k}k}{(2k+1)16^k}\sum_{j=0}^k\f1{(2j+1)^3}\eq\f7{180}\l(\f{-1}p\r)pB_{p-5}\pmod{p^2}.
\end{equation*}
\end{conjecture}

\begin{conjecture}\label{Conj5.2} {\rm (i)} We have
\begin{align}\label{5.3}\sum_{k=0}^\infty\f{\bi{2k}k}{(2k+1)8^k}\(\sum_{0\ls j<k}\f{(-1)^j}{2j+1}-\f{(-1)^k}{2k+1}\)=&-\f{\sqrt2}{16}\pi^2,
\\\label{5.4}\sum_{k=0}^\infty\f{\bi{2k}k}{(2k+1)16^k}\(12\sum_{j=0}^k\f{(-1)^j}{(2j+1)^2}-\f{(-1)^k}{(2k+1)^2}\)=&4\pi G,
\\\label{5.5}\sum_{k=0}^\infty\f{\bi{2k}k}{(2k+1)^2(-16)^k}\(5\sum_{j=0}^k\f1{(2j+1)^3}+\f1{(2k+1)^3}\)=&\f{\pi^2}2\zeta(3),
\\\label{5.6}\sum_{k=0}^\infty\f{\bi{2k}k}{(2k+1)16^k}\(24\sum_{0\ls j<k}\f{(-1)^j}{(2j+1)^3}+7\f{(-1)^k}{(2k+1)^3}\)=&\f{\pi^4}{12},
\end{align}
\begin{align}
\label{5.7}\sum_{k=0}^\infty\f{\bi{2k}k}{(2k+1)^2(-16)^k}\(40\sum_{0\ls j<k}\f{(-1)^j}{(2j+1)^3}-7\f{(-1)^k}{(2k+1)^3}\)=&-\f{85\pi^5}{3456},
\\\label{5.8}\sum_{k=0}^\infty\f{\bi{2k}k}{(2k+1)16^k}\(3\sum_{j=0}^k\f1{(2j+1)^4}-\f1{(2k+1)^4}\)=&\f{121\pi^5}{17280},
\\\label{5.9}\sum_{k=0}^\infty\f{\bi{2k}k}{(2k+1)^2(-16)^k}\(5\sum_{0\ls j<k}\f1{(2j+1)^4}+\f 1{(2k+1)^4}\)=&\f{7\pi^6}{7200}.
\end{align}

{\rm (ii)} Let $p>3$ be a prime. Then
\begin{align*}&\sum_{k=0}^{(p-3)/2}\f{\bi{2k}k}{(2k+1)16^k}\(12\sum_{j=0}^k\f{(-1)^j}{(2j+1)^2}-\f{(-1)^k}{(2k+1)^2}\)\eq-\f32\cdot\f{H_{p-1}}{p^2}\pmod{p^2},
\\&\sum_{k=0}^{(p-3)/2}\f{\bi{2k}k}{(2k+1)^2(-16)^k}\(5\sum_{j=0}^k\f1{(2j+1)^3}+\f1{(2k+1)^3}\)\eq\f{B_{p-5}}8\pmod p,
\\&\sum_{k=0}^{(p-3)/2}\f{\bi{2k}k}{(2k+1)16^k}\(24\sum_{0\ls j<k}\f{(-1)^j}{(2j+1)^3}+7\f{(-1)^k}{(2k+1)^3}\)
\\&\qquad\ \eq-\f p5B_{p-5}\pmod{p^2},
\\&\sum_{k=0}^{(p-3)/2}\f{\bi{2k}k}{(2k+1)^2(-16)^k}\(40\sum_{0\ls j<k}\f{(-1)^j}{(2j+1)^3}-7\f{(-1)^k}{(2k+1)^3}\)
\\&\qquad\ \eq\l(\f{-1}p\r)\f 5{32}B_{p-5}\pmod{p^2}.
\end{align*}
Provided $p>7$, we have
\begin{align*}&\sum_{k=0}^{(p-3)/2}\f{\bi{2k}k}{(2k+1)16^k}\(\sum_{j=0}^k\f3{(2j+1)^4}-\f1{(2k+1)^4}\)
\\&\qquad\qquad\eq\f5{192}\l(\f{-1}p\r)\f{H_{p-1}^{(3)}}{p^2}\pmod{p^2}
\end{align*}
and
$$\sum_{k=0}^{(p-3)/2}\f{\bi{2k}k}{(2k+1)^2(-16)^k}\(\sum_{0\ls j<k}\f5{(2j+1)^4}+\f1{(2k+1)^4}\)\eq-\f{pB_{p-7}}{560}\pmod{p^2}.$$
\end{conjecture}

\begin{conjecture}\label{Conj5.3} We have
\begin{equation}\label{5.10}
\begin{aligned}&\sum_{k=0}^\infty\f{\bi{2k}k}{(2k+1)16^k}\(8\sum_{j=0}^k\f{(-1)^j}{(2j+1)^4}+\f{(-1)^k}{(2k+1)^4}\)
\\&\qquad=\f{11}{120}\pi^2\zeta(3)+\f 83\pi \beta(4),
\end{aligned}\end{equation}
\begin{equation}\label{5.11}
\begin{aligned}&\sum_{k=0}^\infty\f{\bi{2k}k}{(2k+1)16^k}\(\sum_{j=0}^k\f{33}{(2j+1)^5}+\f4{(2k+1)^5}\)
\\&\qquad=\f{35}{288}\pi^3\zeta(3)+\f{1003}{96}\pi\zeta(5),
\end{aligned}\end{equation}
\begin{equation}\label{5.12}
\begin{aligned}&\sum_{k=0}^\infty\f{\bi{2k}k}{(2k+1)^2(-16)^k}\(110\sum_{j=0}^k\f{(-1)^j}{(2j+1)^4}+29\f{(-1)^k}{(2k+1)^4}\)
\\&\ \ \qquad=\f{91}{96}\pi^3\zeta(3)+11\pi^2\beta(4)-\f{301}{192}\pi\zeta(5),
\end{aligned}\end{equation}
\begin{equation}\label{5.13}
\begin{aligned}&\sum_{k=0}^\infty\f{\bi{2k}k}{(2k+1)^316^k}\(72\sum_{j=0}^k\f{(-1)^j}{(2j+1)^2}-\f{(-1)^k}{(2k+1)^2}\)
\\&\ \ \qquad=\f{7}{3}\pi^3G+\f{17}{40}\pi^2\zeta(3),
\end{aligned}\end{equation}
\begin{equation}\label{5.14}
\begin{aligned}&\sum_{k=0}^\infty\f{\bi{2k}k}{(2k+1)^316^k}\(\sum_{j=0}^k\f{33}{(2j+1)^3}+\f8{(2k+1)^3}\)
\\&\qquad=\f{245}{216}\pi^3\zeta(3)-\f{49}{144}\pi\zeta(5).
\end{aligned}\end{equation}
\end{conjecture}

\end{document}